\documentstyle[12pt,euscript,leqno]{amsart}
\newtheorem{thm}{Theorem}[section]
\newtheorem{lem}[thm]{Lemma}
\newtheorem{cor}[thm]{Corollary}

\numberwithin{equation}{section}

\newcommand{\Om}{\Omega}

\newcommand{\RR}{{\mathbb{R}}}
\newcommand{\Ga}{\Gamma}

\renewcommand{\phi}{\varphi}

\newcommand{\UU}{{\mathcal U}}

\newcommand{\DD}{{\mathcal D}}

\begin{document}

\title[PRESCRIBED JACOBIAN EQUATIONS]
{ON THE LOCAL THEORY OF PRESCRIBED JACOBIAN EQUATIONS REVISITED}
\author{Neil S. Trudinger}
\thanks{Research supported by Australian
Research Council Grants DP170100929, DP180100431.}
\address{${}^\dagger$Mathematical Sciences Institute,
Australian National University, Canberra, ACT 0200, Australia.}
\email{neil.trudinger@@anu.edu.au}

\maketitle

\vspace{0.3cm}

\noindent { Abstract.}
\vspace{0.3cm}

In this paper we revisit our previous study of the local theory of prescribed Jacobian equations associated with generating functions, which are extensions of cost functions in the theory of optimal transportation. In particular, as foreshadowed in the earlier work,  we provide details pertaining to the relaxation of a monotonicity condition in the underlying convexity theory and the consequent classical regularity. Taking advantage of recent work of Kitagawa and Guillen, we also extend our classical regularity theory to the weak form A3w of the critical matrix convexity conditions.

\vspace{1cm}

%
%
\section{Introduction}
In this paper we revisit our previous study \cite {T5} of the local theory of prescribed Jacobian equations associated with generating functions, which are extensions of cost functions in the theory of optimal transportation. In particular we elaborate further our remark there pertaining to relaxing  the monotonicity condition on the matrix function $A$ in our convexity theory, thereby enabling the use  of duality properties in the ensuing convexity and local regularity theory. 

We begin by describing the class of equations under consideration, which we called generated prescribed Jacobian equations and is now typically abbreviated to just \emph{generated Jacobian equations},  (GJEs). Let $\Om$ be a domain in Euclidean $n$-space, $\mathbb{R}^{\it n}$, and
$Y$ a $C^1$ mapping from $\Om \times\mathbb{R} \times\RR^{\it n}$ into $\mathbb{R}^{\it n}$. The
\emph{prescribed Jacobian equation}, (PJE), is a partial differential equation of
the form,
%
%
\begin{equation} \label{def:PJE}
 \det DY (\cdot,u,Du)  =  \psi(\cdot,u,Du),
\end{equation}

\noindent where $\psi$ is a given scalar function on
$\Om\times\RR\times{\RR}^n$ and $Du$ denotes the gradient of the
function $u:\Om \to\RR$.

Denoting points in $\Omega\times\RR\times\RR^n$ by $(x,z,p)$, we always assume that the matrix $Y_p$ is invertible,
that is $\det Y_p \neq 0$, so that we may write \eqref{def:PJE} as a general
\emph{Monge-Amp\`ere type equation}, (MATE),
%
%
\begin{equation} \label{MATE}
\det [D^2u - A(\cdot,u,Du) ] =  B(\cdot,u,Du),
\end{equation}

\noindent where
%
%
\begin{equation} \label{def:AB}
A = - Y^{-1}_p(Y_x+Y_z\otimes p),
\quad B = (\det Y_p)^{-1}\psi.
\end{equation}

\noindent A function $u\in C^2(\Om)$ is degenerate elliptic,
(elliptic), for equation \eqref{MATE}, 
whenever
%
%
\begin{equation} \label{def:elliptic}
D^2u - A(\cdot,u,Du) \ge 0, \quad (>0),
\end{equation}
\noindent in $\Om$. If $u$ is an elliptic solution of \eqref{MATE}, then the function $B(\cdot,u,Du)$
is positive. Accordingly we  assume throughout that $B$ is at least non-negative
in $\Om\times\RR\times\RR^n$, that is $\psi$ and $\det Y_p$ have the same
sign.

 The \emph{second boundary
value problem} for the prescribed Jacobian equation is  to prescribe the image,
%
%
\begin{equation} \label{def:bvp}
Tu(\Om): = Y (\cdot,u,Du)(\Om) = \Om^*,
\end{equation}

\noindent where $\Om^*$ is another given domain in
 $\RR^n$. When
$\psi$ is separable, in the sense that
%
%
\begin{equation} \label{def:separable}
|\psi(x,z,p)| = f(x)/f^*\circ Y(x,z,p),
\end{equation}
\noindent for positive $f, f^*\in L^1(\Om)$, $L^1(\Om^*)$ respectively, then a
necessary condition for the existence of an elliptic solution,
for which the mapping $T$ is a diffeomorphism, to
the second boundary value problem \eqref{def:PJE}, \eqref{def:bvp} is the \emph{mass
balance} condition,
%
%
\begin{equation} \label{def:massbalance}
\int_\Om f=\int_{\Om^*} f^*.
\end{equation}
\vspace{0.3cm}

We note that $Y$ need only be defined on the one jet of $u$, $J_1 = J_1[u] (\Om)= (\cdot,u,Du)(\Om)$
in order to formulate (1.1) and (1.5) and typically we will only have $Y$ defined on an open set
$\UU \subset \RR^n \times \RR\times\RR^n$ with the resultant Monge-Amp\`ere type equation (1.2), 
accompanied by a constraint, $J_1[u] \subset\UU$. 

We will also change our notation slightly from \cite{T5} 
and let $g\in C^2(\Gamma)$ denote a \emph{generating function}, where  $\Gamma$ is a domain in $\RR^n\times\RR^n\times\RR$ whose projections, $$I(x,y) = \{z\in\RR |\  (x,y,z)\in\Gamma\},$$
are open intervals. For convenience, we also denote the following projections,
$$\Gamma_x = \{(y,z) \in \RR^n\times\RR |\  (x,y,z)\in\Gamma \},\quad \Gamma_{y,z} = \{x\in \RR^n |\  (x,y,z)\in\Gamma \}.$$
Note that these projections may be empty  for some values of $x,y$ and $z$.

Denoting points in $\RR^n\times\RR^n\times\RR$, by $(x,y,z)$, we assume that $g_z\ne 0$ in $\Gamma$, 
together with the following two conditions which extend the corresponding
conditions in the optimal transportation case, \cite{MTW}:

\vspace{0.2cm}

\begin{itemize}
\item[{\bf A1}:]
The mapping $(g_x,g)(x,\cdot,\cdot) $ is one-to-one
in $\Gamma_x$, for each $x\in \RR^n$.
\vspace{0.2cm}

\item[{\bf A2}:]

$\det E\neq 0$ in $\Gamma$, where $E$ is the $n\times n$ matrix given by

%
%
\begin{equation} \label{def:E}
E =  [E_{i,j}] = g_{x,y} - (g_z)^{-1}g_{x,z}\otimes g_y.
\end{equation}
\end{itemize}
\vspace{0.3cm}
From A1 and A2, the vector field $Y$, together with the \emph{dual} function $Z$,
are generated by $g$ through the equations,

%
%
\begin{equation} \label{def:YZ}
g_x( x, Y, Z) = p,\quad g(x,Y,Z) =u.
\end{equation}
\vspace{0.3cm}

\noindent 
The Jacobian determinant of the mapping $(y,z) \to  (g_x,g) (x,y,z)$
is  $g_z \det E, \neq 0$ by A2, so that $Y $and $Z$ are accordingly $C^1$ smooth.
Also by differentiating \eqref {def:YZ}, with respect to $p$, we obtain $Y_p = E^{-1}$.
Using \eqref{def:AB} or differentiating \eqref{def:YZ} for $p = Du$, with respect to $x$, 
we obtain that
the corresponding prescribed Jacobian equation \eqref{def:PJE}  is a Monge-Amp\`ere equation of the
form \eqref {MATE}  with

%
%

\begin{align}\label{def:GPJE}
A(x, u, p)  &= g_{xx}[x,Y(x,u, p),Z(x,u,p)] , \\
 B(x, u, p) &=\det E(x,Y,Z)\psi (x, u, p) \notag
\end{align}

\noindent and is well defined in domains $\Om$ for $J_1 = J_1[u] (\Om)\subset\UU$,
where
%
%
\begin{equation}
 \UU =  \big \{(x,u,p)\in\RR^n\times\RR\times\RR^n \mid u = g(x,y,z), p = g_x(x,y,z), (x,y,z) \in \Ga\big\}.
 \end{equation}
\
 
 Following \cite{T5} we also have the dual condition to A1:
 \vspace{0.2cm}
\begin{itemize}
\item[{\bf A1*}:]

The mapping $Q: = -g_y/g_z(\cdot,y,z)$ is one-to-one in $\Gamma_{y,z} $, for all $(y, z)\in\RR^n\times\RR$.

\end{itemize}
\vspace{0.2cm}

 Condition A1* arises through the notion of duality introduced in \cite{T5}, where 
 the dual generating function $g^*$  is defined by 
%
%
\begin{equation}\label {def:g^*}
g[x,y,g^*(x,y,u)] = u.
\end{equation}
\noindent Clearly $g^*$ is well defined on the dual set,

$$\Gamma^* = \big\{(x,y,u) \in \RR^n\times\RR^n\times\RR \mid u\in J(x,y)\big\},$$

\noindent where $J(x,y) = g(x,y,\cdot)I(x,y)$, and $g^*_y (x,y,u) = Q(x,y,z)$ for $u=g(x,y,z)$ so that condition
A1* may be equivalently expressed as the mapping $g^*_y$ is one-to-one in $x,u$ for all $(x, y,
u)\in \Ga^*$. Furthermore the Jacobian matrix of the mapping $x\to Q(x,y,z)$ is $-E^t/g_z$, 
where $E^t$ denotes the transpose of $E$,  so its determinant is
automatically non-zero when condition A2 holds. From condition A1*, we then infer the existence of a $C^1$ dual mapping $X$,
 defined uniquely by
 %
 %
 \begin{equation} \label{def:X}
 Q\Big(X(y,z,q), y, z\Big) = q
 \end{equation}
 \noindent for all $q\in Q(\cdot,y,z)(\Gamma_{y,z})$. Note also that by setting 
$$P(x,y,u) = g_x\big(x,y,g^*(x,y,u)\big),$$
\noindent we may express condition A1 in the same form as  A1*, namely the mapping $P$ is one-to-one in $y$, for all
$(x,u)$ such that ($x,y,u)\in \Ga^*$.

In the special case of \emph{optimal transportation}, we have
%
%

\begin{align}\label{def:OTG}
&g(x,y,z) = -c(x,y) - z,\quad \Ga = \DD\times\RR\quad g_z = -1,\quad I = I(x,y) = J(x,y) = \RR,\\
& E =- c_{x,y}, \quad g^*(x,y,u) =- c(x,y) -u, \notag
\end{align}
\noindent where $\DD$ is a domain in  $\RR^n\times\RR^n$ and $c\in C^2(\DD)$ is a \emph{cost function}, satisfying conditions A1 and A2 in \cite{MTW}. The essential difference here is that $Y$ and $A $ are independent of $u$ so that our arguments here and in \cite{T5} are primarily concerned with handling such a dependence. 

As in \cite{T5} we will assume throughout that $g$  has been normalised so that $g_z <0$ in accordance with \eqref{def:OTG}.

Our next conditions extend the conditions A3 and A3w introduced for optimal transportation in \cite{MTW,T1,TW3}
and are expressed in terms of the matrix function $A$ in \eqref{MATE}, which for the purpose of classical regularity is assumed twice differentiable.

\vspace {0.3cm}

\noindent {\bf A3}\ ({\bf A3w}) \vspace{-0.5cm}
$$A^{kl}_{ij}\xi_i\xi_j\eta_k\eta_l: = (D_{p_kp_l}A_{ij}) \xi_i\xi_j\eta_k\eta_l > (\ge)\ 0,
        $$

\vspace {0.2cm}

\noindent for all $(x,u,p) \in \UU,  \xi,\eta \in \RR^n$ such that $\xi.\eta = 0$.

\vspace {0.3cm}

\noindent Conditions A3w (A3) express a co-dimension one convexity (strict codimension one convexity) 
of the matrix function $A$ with respect to the gradient variable $p$ in the set  $\UU$, which we can generally assume is convex in $p$ for fixed $x$ and $u$. As in \cite{T1}, we may write equivalently that $A$ is \emph{ regular}, 
(\emph{strictly regular}), in $\UU$. It is proved in \cite{T5} that conditions A3 and A3w are invariant under duality, through explicit formulae for
$D_p^2A$ in terms of the generating function $g$ and its derivatives up to order four. This result is extended to non smooth $A$ in \cite{LoT2}, where
$A$ co-dimension one convex (strictly co-dimension one convex) means that the form $A\xi.\xi = A_{ij}\xi_i\xi_j $ is convex, (locally uniformly convex), along line segments in $p$, orthogonal to $\xi$ for all $\xi\in \RR^n$.

In \cite{T5} we also introduced conditions expressing the monotonicity of $A$ with respect to $u$, namely:

\vspace {0.3cm}

\noindent {\bf A4}\ ({\bf A4w}) \vspace{-0.5cm}

$$D_uA_{ij}\xi_i\xi_j > (\ge 0),$$

\vspace {0.2cm} 
\noindent for all $(x,u,p) \in \UU,  \xi \in \RR^n$. 

\vspace {0.2cm}

\noindent Only the weak monotonicity A4w was used in \cite{T5}. 






In the next section, we revisit the corresponding section in \cite{T5} and show that conditions A1, A2, A1* and A3w
suffice for the convexity theory developed there, without using the monotonicity condition A4w. This will entail some upgrading of our conditions on domains, relative to $\Gamma$, but will facilitate better the use of duality in ensuing regularity arguments. In fact we had already 
worked out versions of these out at the time of writing \cite{T5} but omitted them in order to avoid the messier statements which pertained to ensuring the sets where condition A3w is used are contained in $\Gamma$.   

 In Section 3, as foreshadowed in \cite{T5} we revisit our existence and interior regularity theory. We work with a modified version of our gradient control assumption in \cite{T5}. For this and throughout this paper, it will be convenient to fix domains $U$ and $V$ in $\RR^n$ such that $U\times V \times g^*\{ U\times V \times J(U,V)\} \subset \Gamma$, where $J(U,V) = \cap_{U\times V} J$. For domains $\Om$ and $\Om^*$ with $\bar\Om \subset U$ and $\bar\Om^* \subset V$, we will then assume for our existence results:

\vspace {0.2cm}

\noindent {\bf A5}: \ \  There exists  an open interval $J_0 =(m_0,M_0) \subset J(U,V)$, $-\infty\le m_0< M_0\le\infty$ and positive constant $K_0 < (M_0-m_0)/2d$, $d = {\rm diam} (\Om)$, such that 

$$  |g_x(x,y,z)| \le K_0 $$

\vspace {0.2cm} 
\noindent for all $x\in\bar\Om, y\in\bar\Om^*, g(x,y,z) \in J_0$.

\vspace{0.3cm}

Following \cite{T5}, we then have the following classical existence theorem, which improves the corresponding result in Corollary 4.7 there.   For the notions of domain convexity used here, the reader is referred to \cite{T5} or Section 2 of this paper.

\begin{thm} 
Let $\Omega$ and $\Omega^*$ be bounded domains in $ \RR^n$,  and let $g$ be a generating function satisfying A1,A2,A1*, A3 and A5.
 Suppose that $f>0, \in C^{1,1}(\Om), f^* >0, \in C^{1,1}(\Om^*)$, with $f,1/f \in L^\infty(\Om), f^*,1/f^* \in L^\infty(\Om^*)$ and that $f$ and $f^*$ satisfy the mass balance condition  \eqref{def:massbalance}.Then for any $x_0\in\Om$ and $u_0$ satisfying  $m_0 +K_1 < u_0 < M_0 - K_1$ , $K_1 = K_0 \text{diam} (\Omega)$, there exists a $g$-convex, elliptic solution $u\in C^3(\Om)$ of the second boundary value problem \eqref{MATE},\eqref{def:bvp}, satisfying $u(x_0) = u_0$, provided $\Om^*$ is $g^*$-convex with respect to $\Om\times J_1$, where $ J_1 = (u_0 - K_1,u_0+K_1)$, and $\Om$ is $g$-convex with respect to all $y\in \Om^*$ and  $z\in g^*(\cdot,y,J_1) (\Om)$.
 \end{thm}
 
Theorem 1.1 is an immediate consequence of the local regularity result Theorem 3.2, which extends Theorem 4.6 in \cite{T5}. Taking account of recent developments, notably the strict convexity result in \cite{GK}, we can now extend Theorem 1.1 to A3w, by extending the regularity argument in \cite{LT1} for the optimal transportation case and moreover by \cite {Ra} the solution is unique. We will also treat this extension in Section 3, (see Theorem 3.4 and Corollaries 3.5, 3.6),  together with the necessary local Pogorelov estimate for its proof, Lemma 3.3. These results have also been presented by us in recent lectures, at Peking University in 2019 and Okinawa Institute of Science and Technology in early 2020.
 
Finally in Section 4, we revisit again our convexity theory, providing an extension of Theorem 3.2 to non-smooth densities and a variant of our key convexity property Lemma 2.2, which does not need duality for its proof. 

We conclude this introduction by noting that our introduction of the concept of generating function in \cite{T4,T5} was to provide a framework for extending the theory of optimal transportation to embrace near field geometric optics, where the associated ray mappings depended also on the position of a reflecting or refracting surface as well as its gradient. Particular motivation came from the point source reflection regularity theory in \cite{KW} which for graph targets is modelled by the generating function in equation (4.15) in \cite{T5}. Note that it is $-A$ in \cite{T5} (4.17) which satisfies A4 for $\tau<0$ so the local regularity theory in this case is covered here. The reader is also referred to the papers \cite{GK, JT1} for further examples of generating functions in optics which fit our theory here.

\vspace{1cm}

%
%

\section{Convexity theory}


We begin by repeating  the definitions in \cite{T5}. We consider bounded domains $\Om$ and $\Om^*\subset \RR^n$ and
 a generating function $g$, satisfying conditions A1 and A2 on $\Ga$. For $x_0,y_0\in \RR^n$, 
 we also denote
$$I(\Om,y_0)  = \cap_\Om I(\cdot,y_0), \quad J(x_0,\Om^*) =  \cap_{\Om^*}J(x_0,\cdot).$$

\vspace{0.2cm}

 A function $u\in C^0(\Om)$ is called $g$-\emph{convex} in  $\Om$,
 if for each $x_0\in\Om$, there exists $y_0\in \RR^n$ and
 $z_0\in I(\Om,y_0)$
 such that  
 %
%
\begin{align}\label{def:g-convex}
u(x_0) &= g(x_0,y_0,z_0), \\
 u(x) &\ge g(x,y_0,z_0)\notag
  \end{align}
\noindent for all $x\in\Om$. If $u(x)> g(x, y_0,z_0)$ for all $x\ne x_0$, then $u$ is called \emph{strictly} $g$-\emph{convex}. If $u$ is differentiable at $x_0$, then $y_0 = Tu (x_0): = 
Y(x_0,u(x_0),Du(x_0))$,
while if $u$ is twice differentiable at $x_0$, then
\begin{equation}\label{2.2}
 D^2u(x_0) \ge  g_{xx}(x_0,y_0,z_0) = A(\cdot,u,Du)(x_0)
\end{equation}
\noindent that is, $u$ is degenerate elliptic for equation \eqref{def:elliptic}  at $x_0$. If $u\in C^2(\Om)$,
we call $u$ \emph{locally} $g$-\emph{convex} in $\Om$ if $J_1[u](\Omega)\subset\mathcal U$ and \eqref{2.2}
 holds for all $x_0\in\Om$.
We will also refer to functions of the form $g(\cdot,y_0,z_0)$ as \emph{$g$-affine} and as
a \emph{$g$-support} at $x_0$ if \eqref{def:g-convex} is satisfied. Note also that the $g$-convexity of a function
 $u$ in $\Omega$ implies its local semi-convexity. When the inequality in \eqref{2.2} is strict, that is $u$ is elliptic
for equation \eqref{def:elliptic}, then we call  $u$ \emph{locally uniformly} $g$-\emph{convex}. It follows readily  that a locally $g$-convex $C^2$ function $u$ will also be locally $g$-convex in the sense that $u\ge g(\cdot, y_0,z_0)$ in some neighbourhood of $x_0$, while a locally uniformly $g$-convex function is strictly $g$-convex in some neighbourhood of  $x_0$.
 
\vspace{0.2cm}

 The domain $\Om$ is $g$-\emph{convex} with respect to $y_0\in \RR^n$, $ z_0\in I(\Om,y_0)$ 
 if the image $Q_0(\Om): = -g_y/g_z (\cdot,y_0,z_0)(\Om)$ is convex in $\RR^n$.  
  
 \vspace{0.2cm}

The domain $\Om^*$ is $g^*$-\emph{convex}
with respect to $x_0 \in \RR^n$, $u_0 \in J(x_0,\Om^*)$, if 
the image $P_0(\Om^*): = P(x_0,\cdot,u_0)(\Om^*) =g_x [x_0,\cdot,g^*(x_0,\cdot,u_0)](\Om^*)$ is convex in $\RR^n$.
 
\vspace{0.2cm} 

We may also consider a corresponding notion of domain convexity when $u$ is fixed which agrees with that associated with the vector field $Y$ in \cite{T1}. Namely, the domain $\Om$ is $Y$-\emph{convex} with respect to 
$y_0\in \RR^n$, $ u_0\in J(\Om,y_0)$ if 
the image $Q^*_0(\Om): = Q^*(\Om, y_0, u_0) = -g_y/g_z [\cdot,y_0,g^*(\cdot,y_0,u_0)](\Om)$ is convex in $\RR^n$.
 It follows then that
$\Om$ is $g$-convex with respect to $y_0\in \RR^n$, $ z_0\in I(\Om,y_0)$ if $\Om$ is $Y$-convex with respect to $y_0$ and $u_0 = g(x,y_ 0, z_0)$ for every $x\in \Omega$. We remind the reader that our definition of g*-convexity is already a special case of the notion of $Y^*$-convexity in \cite{T1} since $P_0(\Om^*) = \{p\in\RR^n\mid Y(x_0,u_0,p)\in\Om^*\}$.

 \vspace{0.2cm}

It will also be convenient to introduce a more  general "sub convexity" notion  as follows. 
The domain $\Om$ is \emph{sub $g$-convex} with respect to $y_0\in \RR^n$, $z_0\in I (\Om,y_0)$   if
the convex hull of  $Q_0(\Om)\subset Q(\Gamma)$ and the domain $\Om^*$ is  \emph{sub $g^*$-convex}
 with respect to  $x_0 \in \RR^n$,  $u_0 \in J(x_0,\Om^*)$  if the convex hull of 
 $P_0(\Om^*)\subset P(\Gamma^*)$. Analogously the domain $\Om$ is \emph{sub $Y$-convex} with respect to 
 $y_0\in \RR^n$, $ u_0\in J(\Om,y_0)$ if 
the image $Q^*_0(\Om)\subset Q^*(\Gamma^*)$.  
 \vspace{0.2cm}

We also define the domain $ \Om^*$ to be $g^*$-convex (sub $g^*$-convex) with respect  to a function $u \in C^0(\Om)$ if $\Om^*$ is $g^*$-convex (sub $g^*$-convex) with respect to each point on the graph of $u$. 

\vspace{0.2cm}

Note that the above definitions also can be applied to general sets, in place of the domains $\Om$ and  $\Om^*$.

\vspace{0.2cm}
 
 Next we define the relevant notions of normal mapping and section.
 
 \vspace{0.2cm}

Let $u\in C^0(\Om)$ be $g$-convex in $\Om$. We define the \emph{$g$-normal} mapping
of $u$ at $x_0\in\Om$ to be the set:

$$Tu(x_0) = \big\{ y_0\in \RR^n \mid 
  \Omega \subset \Gamma_{y_0,z_0} \text{ and } \quad u(x) \ge g(x,y_0,z_0)\text{ for all } x\in\Om\big\},$$
 
\noindent where $z_0 = g^*(x_0,y_0,u_0), u_0 =u(x_0)$. Clearly $Tu$ agrees with our previous terminology when $u$ is differentiable. In the non differentiable case we at least have the inclusion,
  \begin{equation}\label{(2.5)}
 Tu (x_0) \subseteq \Sigma_0: = Y(x_0,u(x_0),\partial u(x_0)),
\end{equation}
where $\partial u$ denotes the subdifferential of $u$, provided the extended one jet, 
$J_1[u](x_0) = [x_0,u(x_0), \partial u(x_0)]\subset \mathcal U$. Moreover from the semi-convexity of $u$, it follows that
$\partial u(x_0)$ is the convex hull of $P_0(Tu(x_0))$ and dist$\{Tu(x),Tu(x_0)\} \rightarrow 0$ as $x\rightarrow x_0$.

 \vspace{0.2cm}
   
Next if $g_0 = g(\cdot,y_0,z_0)$ is a $g$-affine function on $\Omega$, we define the \emph{$g$-section} $S$ of  a $g$-convex function $u$ with respect to $g_0$ by
 
$$ S = S(u,g_0) = S(u,y_0,z_0) = \big\{ x\in\Om \mid u(x) <  g(x,y_0,z_0) \big\} $$

\noindent If $g_0$ is also a $g$-affine support to $u$ at $x_0$, we define the \emph{ contact set} $S_0$ by
 
$$ S_0 = S_0(u,g_0) = S_0(u,y_0,z_0) = \big\{ x\in\Om \mid u(x) = g(x,y_0,z_0) \big\}.$$ 
 
\noindent Note that we have defined sections here differently to \cite{T5}.

 \vspace{0.2cm}
   
We now have the following variant of Lemma 2.1 in \cite{T5}.

  \begin{lem} 
Assume that g satisfies A1, A2, A1* and A3w  and that $u\in C^2(\Om)\cap C^0(\bar\Om)$
is locally $g$-convex in $\Om$ and $u(\Om)\subset\subset J(\Om,Tu(\Om))$ with $Tu(\Om)$ sub $g^*$-convex
with respect to $u$.  Then if $\Om$ is $g$-convex with respect
 to $(y,z)$  for all $y \in Tu(\Om), z \in g^*(\cdot,y,u)(\Om)$, it follows that  $u$ is $g$-convex in $\Om$.
\end{lem}

\noindent {\bf Remark 2.1.}  
More specifically, we have for any $x_0\in\Om, y_0 = Tu(x_0), z_0 =  g^*(x_0,y_0,u(x_0))$, the
$g$-affine function $g(\cdot ,y_0,z_0)$ is a $g$-support, provided  $g(\cdot, y_0, z_1) \ge u$ in $\Om$
for some $ z_1 < z_0, \in I (\Om,y_0)$ and  $\Om$ is $g$-convex
with respect to $(y_0,z)$, for all $z\in ( z_1,z_0)$. 
More generally we can weaken the assumption that $u(\Om)\subset\subset J(\Om,Tu(\Om))$ to just requiring $(\cdot, Tu(\Om), u)(\bar\Om) \subset \Gamma^*$.

\vspace {0.2cm}

In order to prove Lemma 2.1 and the ensuing results concerning $g$-normal mappings and sections, we first recall
a fundamental inequality from \cite{T5}, for which we assume the generating function $g$ satisfies conditions A1, A2, A1*
and A3w. Let $u\in C^2(\Omega)$ be locally g-convex in $\Omega$ and $g_0 = g(\cdot,y_0,z_0)$
be a $g$-affine function defined on $\Omega$, where the domain $\Omega$ is assumed to be $g$-convex with respect to 
$(y_0,z_0)$. Defining the height function $h= u-g_0$ and making the coordinate transformation $x\to q = Q(x,y_0,z_0)$, we then have,
following the computation in  (2.9) and (2.10) in \cite{T5}, (without using condition A4w),  the differential inequality,

\begin{equation}\label{diff ineq}
D_{q_\xi q_\xi}h\ge -K|D_{q_\xi}h| -K_0|h|,
\end{equation}
for any unit vector $\xi$, at any point $\hat x \in \Om^\prime \subset\subset \Om$, for which  $(\cdot, [g_0, u] \cup [u,g_0], Dg_0)(\hat x) \subset \mathcal U$ and the set $Tu(\Omega)\cup\{y_0\}$,  is sub $g^*$-convex with respect to $(\cdot, g_0)(\hat x)$, where $K$ and $K_0$ are constants depending on $g, g_0,\Omega, \Om^\prime$ and $J_1[u]$. If additionally condition A4w holds, the differential inequality \eqref{diff ineq} holds with $K_0=0$, for $h(\hat x) \ge 0$, as in inequality (2.10) in \cite{T5}, while if A3w is replaced by the strict condition A3 or $u$ is locally uniformly $g$-convex at $\hat x$, we have  strict inequality in \eqref{diff ineq},  and our proofs here can be simplified by only using the simpler inequalities

\begin{equation}\label{simpler diff ineq}
D_{q_\xi q_\xi}h > 0
\end{equation}
whenever $h(\hat x) = D_{q_\xi}h(\hat x) = 0$.

 \vspace{0.2cm}

Lemma 2.1 now follows by a modification of the proof of the corresponding lemma in \cite{T5}, which we indicate after
formulating the extensions to normal mappings and sections.  The crucial convexity property for our regularity considerations is the characterisation of the $g$-normal mapping in terms of the sub differential through equality in the inclusion \eqref{(2.5)}. In our earlier versions of this paper going back to 2014, we had various hypotheses for this result depending on whether we raise or lower $g$-affine functions or use duality in the proofs.  The following version, depending on duality, will be convenient for our purposes here.

\begin{lem}
Assume $g$ satisfies A1, A2, A1* and A3w  and suppose $u\in C^0(\Om)$ is $g$-convex in $\Om$, with  $\Omega$ sub $g$-convex  with respect to all $y \in \Sigma_0$, $z=z_y$ for some $x_0\in \Om$. Then we have $Tu(x_0) = \Sigma_0$.
\end{lem}

We will deduce Lemma 2.2 from a convexity result for $g$-sections which extends Lemma 2.3 in \cite{T5}.

 \begin{lem}
Assume $g$ satisfies A1,A2,A1* and A3w  and suppose $u\in C^0(\Om)$ is $g$-convex in $\Om$.  Assume also:
\vspace{2mm}

(i) $\Omega$ is $g$-convex with respect to some $y_0\in \RR^n$, $ z_0\in I(\Om,y_0)$;
\vspace{2mm}

(ii) $Tu(\Omega)\cup \{y_0\}$  is sub $g^*$-convex with respect to 
$g_0: = g(\cdot, y_0,z_0)$ in $\Omega$.
\vspace{2mm}

\noindent Then the section $S = S(u,g_0)$ is also $g$-convex with respect with respect to $(y_0, z_0)$, while if $g_0$ is a $g$-support to $u$, the contact set $S_0 = S_0(u,g_0)$ is $g$-convex with respect with respect to $(y_0, z_0)$. 

\end{lem}

We note that Lemmas 2.2 and 2.3 are direct extensions of the corresponding lemmas in \cite{T5} under condition A4w but Lemma 2.1 needs stronger domain convexity conditions in its hypotheses. We will consider further variants of our convexity results in Section 4, including a  version of Lemma 2.2 which does not need duality in its proof. We also remark here that Lemmas 2.2 and 2.3 are also extended to $C^2$ generating functions in \cite{LoT2}.

\vspace{0.5cm}

\noindent {\bf Proofs: } We indicate here the necessary modifications of the corresponding proofs of Lemmas 2.1 and 2.3 in \cite{T5} which follow by using the more general differential inequality \eqref{diff ineq} in conjunction with decreasing or increasing the appropriate $g$-affine functions.  Letting $u$ be locally $g$-convex in $\Om$, for a fixed point $x_0\in\Om$, $y_0=Tu(x_0)$, $z_0 = g^*(x_0,y_0,u_0)$, we modify the height function $h$ in equation (2.8) in \cite{T5} by setting for $\delta\ge0$,
 \begin{equation}\label{height}
 h(x) = h_\delta(x) =u(x) - g(x,y_0,z_0 - \delta).
 \end{equation}
 As in \cite{T5}, the function $g_0 = g(\cdot, y_0,z_0)$ is a local support near $x_0$, that is $h_0 \ge 0$ near $x_0$. 
 If $h_0(x) < 0$ at some point $x\in \Omega$, from the hypotheses of Lemma 2.1, there exists $\delta \ge 0$  such that $h_\delta$ attains a zero maximum along the closed $g$-segment joining $x_0$ and $ x$, with respect to $y_0,z_0 -\delta$, (which will lie in some subdomain $\Om^\prime \subset\subset \Om$) and  $ h_\delta(x) < 0$. Setting $q_0 = Q(x_0,y_0,z_0 - \delta)$, $q_t = tq + (1-t)q_0$, $x_t = X(q_t,y_0,z_0 -\delta)$, $0\le t\le1$ and defining  $f(t) = h_\delta(x_t)$, it follows that $f$ attains a zero maximum at some point $\hat t \in (0,1)$ for $\delta > 0$ and  $\hat t=0$ for $\delta = 0$, so that also $f^\prime(\hat t) = 0$ in both cases.  From the differential inequality \eqref{diff ineq} and our sub $g$-convexity hypothesis, we have the corresponding differential inequality for $f$,

\begin {equation}  \label {f diff ineq}
f^{\prime\prime} \ge -K|f^\prime| - K_0 |f|
\end{equation}                                                             

\noindent holding, at least in some neighbourhood of the set where $f$ vanishes, which is a contradiction. The latter assertion is easily seen as \eqref{f diff ineq} clearly implies a uniform bound from above for $v^{\prime\prime}$ where $v = \log (\epsilon - f)$, $\epsilon > 0$, and is a one dimensional version of the strong maximum principle.  Accordingly $g_0$ must be a $g$-support and Lemma 2.1 is proved.

 \vspace{0.2cm}

Next we show that Lemma 2.3  also can be proved by a corresponding modification of the proof of Lemma 2.3 in \cite{T5}. By modifying $\Omega$ we may assume if necessary that $u$ and $g_0$ extend to a neighbourhood of $\bar\Omega$.   If $S$ is not $g$-convex, with respect to $y_0, z_0$, there must be a $g$-segment in $\bar\Omega$ joining two points in $S$ containing a  point $x_1\in \partial S$. Setting $ u_1 = u(x_1), y_1\in Tu(x_1), z_1 = g^*(x_1,y_1,u_1)$, the inequality 
\begin{equation}
g(x,y_1,z_1) < g(x,y_0,z_0 ),
\end{equation}

\noindent holds for all $x\in S$.  Now replacing $h$  in \eqref{height} by
\begin{equation}
h(x) =  h_\delta(x)  = g(x,y_1,z_1 + \delta) - g(x,y_0,z_0), 
\end{equation}

\noindent where $g$ is extended so that $g_\delta(x) :=g (x,y_1, z_1 + \delta) = -\infty$ where $(x,y_1, z_1 + \delta) \notin \Gamma$, we again obtain  a contradiction with the differential inequality \eqref{diff ineq} for some $\delta \ge 0 $. Note that in this case we are lowering the support  $g(\cdot, y_1,z_1)$ by increasing $z_1$ and also using condition (ii) in the sense that $(x,y_1,z_1 + \delta) \in \Gamma$ whenever $g_\delta = g_0$. 
 
  Finally, the $g$-convexity of the contact set $S_0$ follows by replacing $S$ by $S_0$ in the above proof, (or alternatively using $S_0 = \cap_{\sigma>0} S_\sigma$, where $S_\sigma = S (u, y_0, z_0 - \sigma)$), and hence Lemma 2.3 is proved.
  
  \vspace{0.2cm}
 \noindent {\bf Remark 2.2.}
In the case when $u$ is also a $g$-affine function  $g_1 = g(\cdot, y_1,z_1)$, we can improve Lemma 2.3 so that $\Omega$ can be replaced by a closed set, with possibly empty interior, and the section $S$ replaced by the closed section $\tilde S = \tilde S( g_1,g_0) =   \{ x\in\Om \mid g_1 \le  g_0 \} $, thereby inferring that $\tilde S$ is $g$-convex with respect to $g_0$. In this case we can take $x_1 \in \Omega - \tilde S$ in the proof since $Tg_1 = y_1$.

\vspace{0.2cm}

Now we can prove Lemma 2.2.  Note that we can infer a version of Lemma 2.2  from Lemma 2.3 through the  $g$-transform $v$, given by
%
%
\begin{equation} \label{def:funcv}
v(y) = u^*_g(y) = \text{sup}_{\Om} \ g^*(\cdot,y,u)
\end{equation}

\noindent as the $g^*$-convexity of $Tu(x_0)$ is equivalent to that of the contact set $S^*_0[v, x_0, u_0]$ of $v$. Here we will proceed somewhat differently with the technicalities.  We begin by  taking two $g$ affine supports to $u$ at $x_0$, $g_0 = g(\cdot, y_0,z_0), g_1 = g(\cdot, y_1,z_1)$ and fixing a point $x = x_1$ in $\Omega$, where $g_1 \le g_0$. Letting $y = y_\theta = Y(x_0, u_0, p_\theta)$, for $p_\theta = \theta Dg_0(x_0) + (1-\theta) Dg_1(x_0)$, $\theta\in [0,1]$ denote a point on the closed $g^*$-segment $I_0$,  with respect to $(x_0, u_0)$, joining $y_0$ and $y_1$,  we now define  $u_y = g (x_1,y_0,z_y)$.  We then obtain, from Remark 2.2 and the $g^*$-convexity of $ I_0$, that $g^*(x_1,y_\theta,  u_\theta) \le g^*(x_0, y, u_0) = z_y$  and hence  $g(x_1, y, z_y) \le g(x_1,y_0,z_0)$, for all $y\in I_0$. 

Consequently we have the following extension of the Loeper maximum principle in optimal transportation \cite{Loe} 
\begin{equation} \label{Loeper}
g(\cdot, y, z_y)  \le \max\{ g_0, g_1\}
\end{equation}

\noindent in $\Om$ for all $y\in I_0$ which, by virtue of the semi-convexity of  $u$ implies that $Tu(x_0)$ is $g^*$- convex with respect to $x_0, u_0$ and hence completes the proof of Lemma 2.2.
 \vspace{0.2cm}

\noindent {\bf Remark 2.3.}
 From the proof of Lemma 2.2 we see that condition (i) can be weakened to just requiring the pair of points $\{x_0, x\}$  is sub $g$-convex  with respect to all $y \in \Sigma_0$, $z=z_y$, for all $x\in \Omega$, that is the $g$-segment joining $x_0$ to any point in $\Omega$ is well defined with respect to all $y \in \Sigma_0$, $z=z_y$. 
 
 \vspace{0.2cm}
 
 We also remark that the technicalities in using the differential inequalities \eqref{diff ineq} in the general A3w case can also be simplified  if we use an approximation of  $g_0$ by a uniformly  $g$-convex function, as in \cite{JT2}, so that we only need the simpler strict inequality $f^{\prime \prime} >0$ from \eqref {simpler diff ineq} whenever $f = f^\prime = 0$,  which requires less smoothness of $g$ for its validity when, additionally to $g\in C^2$ and $A$ regular,  either $u$ is uniformly $g$-convex or $A$ is strictly regular. In fact here the differentiability of $A\xi.\xi$ with respect to $p$ in directions orthogonal to $\xi$ would suffice.  However in using our differential inequality approach, (which for the optimal transportation case goes back to \cite{KM}, with simplified versions in \cite{TW2, V2}), we  would still need some smoothness of the generating function $g$ beyond $C^2$ smoothness.  A substantially different geometric approach to our convexity theory, which has its optimal transportation roots in \cite{TW4}, is presented in \cite{LoT1, LoT2}, where we do not need any derivatives beyond second order. A different analytic approach, based on a weak form of condition A3w corresponding to a sharpening of the quasi-convexity property \eqref{Loeper}, is developed in \cite {GK}.

We will return to the convexity theory in Section 4,  in conjunction with consideration of the strict convexity and continuous differentiability results foreshadowed in \cite{T5}, which are applications of Lemma 2.3.

\vspace{1cm}

%
%

\section{Existence and regularity.}
 
Since the arguments in Section 3 of \cite{T5} are largely independent of condition A4w, they extend readily to the more general case as a consequence of Lemma 2.2.  First we recall from \cite{T5} the definition of generalized solution. For convenience we let $\Omega$ and $\Omega^*$ be bounded domains in $ \RR^n$, whose closures lie in the domains $U$ and $V$ respectively, as introduced in Section 1, and $u\in C^0(\overline\Omega)$ be $g$-convex in  $\Omega$, with $Tu(\Om)\subset\subset V$ and conditions A1,A2, A1* satisfied. Then there is a measure $\mu= \mu[u] = \mu(u,f^*)$ on $\Omega$, for $f^* \ge 0 \in L^1(V)$, such that for any Borel set $E \subset \Omega$, 

%
%
\begin{equation} \label{def:measmu}
\mu (E) = \int _{Tu(E) } f^*
\end{equation}

\noindent which is also weakly continuous with respect to local uniform convergence. A $g$-convex function $u$ on $\Omega$ is now defined to be a generalized solution of the second boundary value problem \eqref{def:bvp} for equation \eqref{def:PJE}, \eqref{def:separable}, under the mass balance condition \eqref{def:massbalance}, if 
%
%
\begin{equation} \label{equ:muf*}
  \mu[u] = \nu_f
\end{equation}

\vspace{0.2cm}

\noindent where $\nu_f = f {\rm d}x$ and $f^*$ is extended to vanish outside $\Omega^*$. We then have the following extension of Theorem 3.1 in \cite{T5}. 

\begin{thm}
Let $\Om$ and $\Om^*$ be domains satisfying  $\bar\Om \subset U$ and $\bar\Om^* \subset V$ and let $g$ be a generating function satisfying A1,A2,A1* and A5.
Suppose that $f$ and $f^*$ are positive densities in $L^1(\Om)$ and $L^1(\Om^*)$ satisfying the mass balance condition \eqref{def:massbalance}.Then for any $x_0\in\Om$ and $u_0 > m_0 +K_1$,
where $K_1 = K_0 \text {diam} (\Omega)$, there exists a generalized solution of \eqref{MATE},\eqref{def:bvp} satisfying $u(x_0) = u_0$.
Furthermore if $\Om^*$ is $g^*$-convex with respect to $u$, then any generalized solution of \eqref{def:bvp} satisfies $Tu(\Om) \subset \bar\Om^*$.
\end{thm}

Using Lemma 2.2 in place of the corresponding Lemma 2.2 in \cite{T5}, we then have the following extension of the local regularity result in Theorem 4.6 in \cite{T5}.

\begin{thm}
Let $u\in C^0(\bar\Om)$ be a generalized solution of \eqref{def:bvp} with positive densities $f \in C^{1,1}(\Om), f^* \in C^{1,1}(\bar\Om^*)$ with $f,1/f \in L^\infty(\Om), f^*,1/f^* \in L^\infty(\Om^*)$ and with generating function $g$ satisfying conditions  A1,A2,A1* and A3.  
Suppose that $u(\Om)\subset\subset J_0$, $\Om^*$ is $g^*$-convex with respect to $u$ and $\Om$ is sub $g$-convex with respect  to the dual function $v = u^*_g$. Then $u\in C^3(\Om)$ and is an elliptic solution of \eqref{MATE}, \eqref{def:GPJE}. Furthermore if $\Om$ is $g$-convex with respect to $v$, then $Tu$ is also a diffeomorphism from $\Om$ to $\Om^*$, with $v$ an elliptic solution of the dual boundary value problem.
\end{thm}

For the explicit form of the dual boundary value problem, the reader is referred to equations (3.5),(3.7) in \cite{T5}. 

 \vspace{0.2cm}

Theorem 1.1 now follows as a consequence of Theorems 3.1 and 3.2, (and approximating $f^*$ near  $\partial\Om^*$ if only locally $C^{1,1}$). From Rankin \cite{Ra}, the solution in Theorem 1.1 is unique. As foreshadowed in \cite{T5}, we will also consider the extension of Theorem 3.2 to $f^* \in C^{1,1}(\Om^*)$ in Section 4, by adapting the  strict convexity argument  in \cite{TW4}, as this also relates to our extension to A3w. 
 \vspace{0.2cm}

In the rest of this section we will treat the extension to A3w which follows for strictly convex generalized solutions from  modification of the Pogorelov estimates in \cite{LTW, TW3}. For this we consider classical elliptic solutions $u$ of the Monge-Amp\`ere type equation \eqref{MATE} in sections $\Om = S(u,g_0)$, with respect to a $g$-affine function $g_0$  on $\bar \Omega$, so that we have a Dirichlet boundary condition $u=g_0$ on $\partial\Om$. Accordingly we  assume $A$ and $B$ are $C^2$ smooth on an  open set $\UU \subset \RR^n \times \RR \times \RR^n$, with $A$ regular and $B$ positive.  We may also assume $\UU$ is convex in $p$, for fixed $x,u$. We then have the following estimate:

\begin{lem}
Let $u\in C^4(\Omega)\cap C^{0,1}(\bar \Omega)$ be an elliptic solution of \eqref{MATE} in $\Omega$ and $u_0 \in C^2(\Omega)\cap C^{0,1} (\bar\Omega)$ a degenerate elliptic solution of the homogeneous equation with $B = 0$, such that  $J_1[u], J_1[u_0] \subset \UU$  and $u=u_0$ on $\partial\Omega$, $u < u_0$ in $\Omega$.  Then there exist positive constants $\beta$, $\delta$ and $C$ depending on $n, \UU, A,B,  u_0$ and  $|u|_1 =\sup_\Omega(|u| + |Du|)$, such that if $d:=\text{diam}(\Omega) < \delta$,

\begin{equation} \label{estimate}
\sup_\Om (u_0 -u)^\beta |D^2u| \le C.
\end{equation}
\end{lem}

For generated Jacobian equations, satisfying A1, A2, A1*, A3w and A4w, Lemma 3.3 corresponds to Theorem 1.2 in \cite{JT1} and does not need the smallness condition on $\Omega$.  We can prove Lemma 3.3 by modification of the global estimate Theorem 3.1 in \cite{TW3}, incorporating an appropriate cut-off function, or by modification of the interior estimate Theorem 2.1 in \cite{LTW}, extending to  $u$ dependence in $A$. Dealing with the general $u$ dependence is the critical issue in the treatment of the cut-off function so we will largely focus on this in the proof. For this it is more convenient for us to begin with the calculations in \cite{TW3}. Accordingly with $u$ and $u_0$ satisfying the hypotheses of Lemma 2.3, we consider, as in \cite{LTW} an auxiliary function 
\begin{equation}
v=v(\cdot, \xi)=\log(w_{ij}\xi_i\xi_j) +\tau |Du|^2+\kappa\phi + \beta \log(u_0 -u)
\end{equation}

\noindent  where $|\xi| =1$,  $\tau$, $\kappa$ and $\beta$ are positive constants to be chosen, $w =D^2u - A(\cdot,u,Du)$ and $\phi\in C^2(\bar\Omega)$ satisfies the global barrier condition,

 \begin{equation}\label{barrier}
      [D_{ij}\phi-D_{p_k}A_{ij}(\cdot,u, Du)D_k\phi]\xi_i\xi_j\geq|\xi|^2
  \end{equation}
  
\noindent Using the smallness condition on $\Omega$ we can fix such a barrier by taking  $\phi = |x-x_1|^2$ for some point $x_1 \in \Omega$.

For the linearized operator, $L$ given by

\begin{equation}\label{linearized operator L}
L=L[u] = w^{ij}[D_{ij}-D_{p_k}A_{ij}(\cdot,u,Du)D_k] - D_{p_k}\log B(\cdot,u,Du)D_k ,
\end{equation}

\noindent with $[w^{ij}]$ denoting the inverse of $[w_{ij}]$,  we can now follow the calculations in \cite{TW3} to obtain, from inequalities (3.11) and (3.12) there, at a maximum point $x_0$ and vector $\xi=e_1$ of $v$ in $\Omega$,

\begin{equation}\label{3.7}
Lv\ge \tau w_{ii} +\kappa w^{ii} - C(\tau+\kappa) + \frac{1}{2w_{11}^2}\sum_{i>1}w^{ii}(D_iw_{11})^2 + \beta L\log(u_0-u),
\end{equation}

\noindent provided $\tau\ge C$ and  $\kappa\ge C \tau$. Here, as is customary, we  use $C$ to denote a constant depending on the same quantities as in the estimate being proved, \eqref{estimate}.  We also note here that a term $\tau w^{ii}$ is missing in the bracketed terms in (3.11), (3.12) in \cite{TW3}, (which is controlled, as in \eqref{3.7}, by taking $\kappa\ge C \tau$). To handle the last term in \eqref{3.7}, we now need to estimate $L\eta$, where $\eta = u_0-u$ by extending analogous estimates in the special cases of optimal transportation in \cite{LTW} and generated Jacobian equations satisfying A4w in \cite{JT1}. We then have in our general case, at $x=x_0\in\Omega$, with similar computation to that underlying our differential inequality\eqref{diff ineq}, using condition A3w,

\begin{equation}\label{3.8}
\begin{array}{rl}
 L \eta \!\!&\!\!\displaystyle = w^{ij}[-w_{ij} + A_{ij}(x,u_0,Du_0) - A_{ij}(x,u,Du) -D_{p_k}A_{ij}(x,u,Du)D_k\eta] \\
  \!\!& \hspace{40ex} \!\!\displaystyle - D_{p_k}(\log B)(x,u,Du)D_k\eta \\
           \!\!&\!\!\displaystyle \ge -C(1+ w^{ii} |\eta| +w^{ii}|D_i\eta|),
\end{array}
\end{equation}
 Substituting  in \eqref{3.7}, we then obtain, with $[w_{ij}]$ assumed diagonal at $x_0$,
\begin{equation}\label{3.9}
Lv\ge \tau w_{ii} +\kappa w^{ii} - C(\tau+\kappa) + \frac{1}{2w_{11}^2}\sum_{i>1}w^{ii}(D_iw_{11})^2 - \frac{C\beta}{\eta} - \frac{C\beta}{\eta^2}w^{ii}|D_i\eta|^2,
\end{equation}
\noindent provided also $\kappa \ge C\beta$. Now we follow \cite{LTW} to control the last term in \eqref{3.9}. Using $Dv(x_0) = 0$, together with $
|D\phi|\le d$, we then estimate, for each  $i=1, \cdots n$,
\begin{equation}
\beta^2 \left(\frac{D_i\eta}{\eta}\right)^2 \le4 \left(\frac{D_iw_{11}}{w_{11}}\right)^2
       +C\tau^2(w_{ii}^2+1)+4\kappa^2 d^2.
\end{equation}
Assuming $\eta w_{11} (x_0) \ge \beta$, we can now estimate the last term in \eqref{3.9},
\begin{equation}
\frac{\beta}{\eta^2}w^{ii}|D_i\eta|^2 \le \frac{C}{\eta} + \frac{4}{\beta w_{11}^2} \sum_{i>1}w^{ii}(D_iw_{11})^2 + C\frac{\tau^2}{\beta} (w_{ii} + w^{ii}) + \frac{4}{\beta}\kappa^2d^2w^{ii}.
\end{equation}
We now conclude the proof of Lemma 3.3 by choosing $\tau \ge C$, $\kappa \ge C\beta$, $\beta \ge C \tau^2$ and $d \le 1/\kappa$.

\vspace{2mm}

 From Lemma 3.3 we now have the following extension of Theorem 3.2 to A3w. 

\begin{thm}
Let $u$ be a strictly $g$-convex generalized solution of \eqref{def:bvp} with positive densities $f \in C^{1,1}(\Om), f^* \in C^{1,1}(\Om^*)$ with $f,1/f \in L^\infty(\Om), f^*,1/f^* \in L^\infty(\Om^*)$ and with generating function $g$ satisfying conditions, A1,A2,A1*and A3w.  
Suppose that  $u(\Om)\subset\subset J_0$, $\Om^*$ is $g^*$-convex with respect to $u$ and $\Om$ is sub $g$-convex with respect to $v = u^*_g$. Then $u\in C^3(\Om)$ and is an elliptic solution of \eqref{MATE}, \eqref{def:GPJE}. Furthermore if  $\Om$ is $g$-convex with respect to $v = u^*_g$, then $Tu$ is also a diffeomorphism from $\Om$ to $\Om^*$, with $v$ an elliptic solution of the dual boundary value problem. 
\end{thm}

To prove Theorem 3.4, we need to adjust the regularity argument in \cite{T5} by using the interior estimate \eqref{estimate} for the solutions $w$ of the approximating Dirichlet problems in Lemma 4.6, in place of the estimate (4.12) in \cite{T5}, when the strong condition A3 is satisfied. To be more specific using Lemma 4.6 in \cite{T5} we construct a sequence $\{\bar u_m\}$ of smooth elliptic solutions of the Dirichlet problem for equation \eqref{MATE} in a small ball of radius $r$, $B_r \subset\Om^\prime \subset\subset \Om$, with boundary condition $\bar u_m = u_m$ on $\partial B_r$, where $\{u_m\}$ is an appropriately chosen sequence of smooth functions converging uniformly to $u$ in $\Om^\prime$, uniformly bounded in $C^1$ and uniformly semi-convex. Here the strict $g$-convexity of $u$ ensures $Tu_m(\Om^\prime) \subset \tilde\Om^*$ for some 
$\tilde\Om^*\subset\subset \Om^*$. Otherwise we need to assume $f^* \in C^{1,1}(\bar\Om^*)$, as done for Theorem 3.2. 
The sequence $\{\bar u_m\}$ will also be uniformly bounded in $C^1$ and moreover with $r$ sufficiently small, $\{\bar u_m\}$  is $g$-convex in $B_r.$  Then using the strict $g$-convexity of $u$ and Lemma 3.3, we can follow the proof of Theorem 3.1 in \cite{LT1} to obtain uniform interior second derivative bounds for $\bar u_m$ in $B_r$. Subsequently, returning to the proof of Theorem 4.6 in \cite{T5} and using Lemmas 4.3 and 4.4 in \cite{T5}, we conclude that $\bar u_m$ converges to $u$ in $B_r$.  From standard elliptic regularity theory \cite{GT}, we then infer $u\in C^3(\Omega)$, is an elliptic solution of \eqref{MATE}, while from the mass balance condition \eqref{def:massbalance}, $Tu (\Om) = \Om^* - E$, for some closed null set $E$. Now letting $T^*v$ denote the $g^*$-normal mapping of $v$ on $\Omega^*$, we then have, (as in the A3 case), from the $g$-convexity of $\Omega$ and smoothness of $u$, that $T^*v(\Om^*) = \Om$. Hence $Tu$ is a diffeomorphism from $\Om$ to $\Om^*$ with $(Tu) ^{-1} = T^*v$ and we conclude from \cite{T5}, Section 3 that $v$ is a $C^3$ elliptic solution of the dual boundary value problem.

\vspace{0.2cm}
In order to apply the strict convexity result of Kitagawa and Guillen \cite{GK}, it will be convenient to strengthen our domain conditions by assuming that the domain $U$ is $g$-convex with respect to all $y \in V$ and $z\in g^*(\cdot, y, J_0)(U)$ and  $V$ is $g^*$-convex with respect to $U\times J_0$. Then we have the following corollaries of Theorem 3.4.

\begin{cor}
Let $u$ be a $g$-convex generalized solution of \eqref{def:bvp} with positive densities $f \in C^{1,1}(\Om), f^* \in C^{1,1}(\Om^*)$ with $f,1/f \in L^\infty(\Om), f,1/f^* \in L^\infty(\Om^*)$, $\bar\Om \subset U$, $\bar\Om^*\subset V$ and with generating function $g$ satisfying conditions, A1,A2,A1*and A3w.  Suppose that $u(\Om) \subset\subset J_0$, together with its $g$-affine supports, and $\Om^*$ is $g^*$-convex with respect to $u$ on $\Om$. Then $u\in C^3(\Om)$ and is an elliptic solution of \eqref{MATE}, \eqref{def:GPJE}. Furthermore if  $\Om$ is $g$-convex with respect to $v = u^*_g$, then $Tu$ is also a diffeomorphism from $\Om$ to $\Om^*$, with $v$ an elliptic solution of the dual boundary value problem.
\end{cor}

Note that from Theorems 2.3 and 2.4 in \cite{GK} we have, in Corollary 3.5, that $u\in C^1(\Omega)$ is strictly $g$-convex in $\Om$, using just the boundedness of the densities $f, f^*$ and their reciprocals.

From Theorem 3.1 and  Corollary 3.5, and taking account of the uniqueness result of Rankin \cite{Ra}, we now have the following extension of Theorem 1.1 to the situation when A3 is weakened to A3w. 

\begin{cor}
Let $\Omega$ and $\Omega^*$ be bounded domains in $\RR^n$, with closures $\bar\Om \subset U$ and $\bar\Om^*\subset V$ and let $g$ be a generating function satisfying A1,A2,A1*, A3 and A5. Suppose that $f>0, \in C^{1,1}(\Om), f^* >0, \in C^{1,1}(\Om^*)$, with $f,1/f \in L^\infty(\Om), f^*,1/f^* \in L^\infty(\Om^*)$ and that $f$ and $f^*$ satisfy the mass balance condition  \eqref{def:massbalance}.Then for any $x_0\in\Om$ and $u_0$ satisfying  $m_0 + 2K_1 < u_0 < M_0 - K_1$, there exists a unique $g$-convex, elliptic solution $u\in C^3(\Om)$ of the second boundary value problem \eqref{MATE},\eqref{def:bvp}, satisfying $u(x_0) = u_0$, provided $\Om^*$ is $g^*$-convex with respect to $\Om\times J_1$, where $ J_1 = (u_0 - K_1,u_0+K_1)$, and $\Om$ is $g$-convex with respect to all $y\in \Om^*$ and  $z\in g^*(\cdot,y,J_1) (\Om)$.
 \end{cor} 
 
 To conclude this section we remark that we can also use Lemma 3.3 to extend the second derivative estimates in \cite{LT1, JT1} and consequently the classical existence theory in\cite{JT2} to generated Jacobian equations under just conditions A1, A1*, A2 and A3w.
 
\vspace{1cm}

\section{More on convexity.}

\noindent {\bf 4.1 Strict convexity}.   In \cite{T5} we also mentioned that the optimal transportation  strict convexity and $C^1$ regularity results in \cite{TW4} extended readily to the generated Jacobian case. The key convexity lemmas for this result are Lemmas 2.3 and its predecessor, Lemma 2.3 in \cite{T5},  on the $g$-convexity of contact sets. 

\begin{lem}
Let $u$ and $g$ satisfy the hypotheses of Lemma 2.3 with $g_0 = g(\cdot,y_0,z_0)$ a $g$-support at $x_0 \in \Om$ and  condition A3w strengthened to A3. Suppose for some ball $B = B_R(x_0) \subset\subset \Omega$, there exists a positive constant $\lambda_0$ such that  
\begin{equation}
|Tu(\omega)| \ge \lambda_0|\omega|,
\end{equation}

\noindent for all open subsets $\omega \subset B$. Then $u$ is strictly $g$-convex at $x_0$. Alternatively suppose that $Tu(\Omega)$ is $g^*$-convex with respect to $u$ on $\Omega$, where  $\Omega$ is sub $g$-convex with respect to $g^*(x,\cdot, u(x))$ on $Tu(\Omega)$ for all $x\in \Omega$, and 
\begin {equation}
|Tu(\omega)| \le \Lambda_0|\omega|
\end{equation}

\noindent for all open subsets $\omega \subset \Omega$ and some positive constant $\Lambda_0$. Then $u \in C^1(\Omega)$.

\end{lem}
To prove Lemma 4.1, we suppose that $u$ is not strictly $g$-convex at $x_0$ so that by Lemma 2.3 there must exist  a $g$-segment $\gamma$, with respect to $y_0,z_0$ joining $x_0$ to another point $x_1\in B$ and lying in $B$. For any sufficiently small radius $\rho$, there then exists a ball $ B_\rho \subset \subset B$ of radius $\rho$ intersecting $\gamma$, with $x_0, x_1 \notin B_\rho$. Now let 
$\tilde u \in C^3(B_\rho)\cap C^0(\bar B_\rho)$ be the unique elliptic solution of the Dirichlet problem,

\begin{equation} \label{4.2}
\det [D^2\tilde u - A(\cdot,\tilde u,D\tilde u) ] =  \frac{\lambda_0}{ 2} \ \text {in} \  B_\rho, \hspace{0.2cm} \tilde u= u \ \text {on}\  \partial  B_\rho, 
\end{equation}

\noindent which, as in the proofs of regularity, is well defined for sufficiently small $\rho$. By our previous comparison arguments, we also have $\tilde u > u \ge g_0$ in $B_\rho$  and $\tilde u = u = g_0$  on $\gamma\cap \partial B_\rho$. This can be seen readily from the finer inequality $\tilde u > \tilde u_0 \ge u$, where  $\tilde u_0$ solves the Dirichlet problem \eqref{4.2} with $\lambda_0/ 2$ replaced by $\lambda_0$. Again Lemma 4.6 in \cite{T5} is also crucial here as in the regularity theory in Section 3. Applying now Lemma 2.3 again in $B_\rho$, with $z_0$ slightly reduced, or more simply the proof of Lemma 2.1, we thus reach a contradiction. The second assertion in Lemma 4.1 on $C^1$ regularity now follows by duality.

Note that in the above argument we have taken advantage of the $g$-convexity of the contact set $S_0$ to simplify the technicalities in \cite{TW4}, where just the connectedness of $S_0$ is used, although the basic approach though reduction to a Dirichlet problem in  small balls is the same. Also, from either the ellipticity of \eqref{4.2} or condition A3,  we only need the strict differential inequality \eqref{simpler diff ineq} in replicating the differential inequality argument in Lemma 2.1. 
 \vspace{0.2cm}

For application to regularity we need to enhance our domain sub  convexity conditions. Namely we will call $\Om$ sub $g$-convex in $U$, with respect to $(y_0, z_0)$, if $z_0 \in I(U,y_0)$ and the convex hull of $Q(\Om,y_0,z_0) \subset Q(U,y_0,z_0)$ and $\Om^*$ sub $g^*$-convex in $V$, with respect to $(x_0, u_0)$, if $u_0 \in J(x_0,V)$ and the convex hull of $P(x_0,\Om^*,u_0) \subset P(x_0,V,u_0)$. Then we can  weaken the $g$-convexity of $\Om$ in Lemma 2.3, in the case of arbitrary contact sets $S_0$, by $\Om$ being sub $g$-convex in a larger domain $U$, with respect to $(y,z)$, for all $g$-supports, $g(\cdot, y,z)$ to $u$, or equivalently with respect to the dual function $v = u_g^*$. In this case we need to strengthen condition (ii) so that $Tu(\Om)$ is sub $g^*$-convex in a larger domain $V$ with respect to $g(\cdot,y,z_y)$ on $\hat\Om_y$,  for all $y\in Tu(\Om), z_y = v(y)$, where $\hat\Om_y$ denotes the $Q(\cdot, y, z_y)$ convex hull of $\Om$. To be more explicit, we define $\hat\Om_y$  as the image under the inverse mapping $Q^{-1}(\cdot,  y,  z_y)$ of the convex hull of $Q(\Om, y, z_y)$. 

Then using duality we have the following extension of Theorem 3.2 to non smooth densities. 

\begin{thm}
Let $u$ be a generalized solution of \eqref{def:bvp} with positive densities $f$ and $f^*$, with $f,1/f \in L^\infty(\Om), f^*,1/f^* \in L^\infty(\Om^*)$ 
and with generating function $g$ satisfying conditions  A1,A2,A1* and A3.  
Suppose that  $\Om\subset\subset U$ is sub $g$-convex in $U$ with respect  to $g^* (x,\cdot, u(x)$ on $Tu(\Om)$, for all $x\in \Om$ and  
$\Om^*\subset\subset V$ is  sub $g^*$-convex in $V$ with respect to $g(\cdot,y,z_y)$ on $\hat\Om_y$ for all $y\in\Om^*, z_y = v(y)$. Then  
$u\in C^1(\bar\Om)$ is strictly $g$-convex in $\Om$ if $\Om^*$ is  $g^*$-convex with respect to $u$ and $Tu$ is a homeomorphism from $\Om$ to 
$\Om^*$ if also  $\Om$ is  $g$-convex  with respect  to $v$.  Furthermore if $f \in C^{1,1}(\Om), f^* \in C^{1,1}(\Om^*)$, then $u \in C^3(\Om)$.
\end{thm}

We can also simplify the  statement of Theorem 4.2, (as well as earlier results),  by replacing $\Gamma^*$ by $U\times V \times J_0$ for sufficiently large domains $U$ and $V$ and interval $J_0$. We note here also that the approach of Loeper \cite{Loe} to local  $C^{1,\alpha}$ regularity in optimal transportation under A3 was used in near field geometric optics in \cite{GuT2015} and has been recently extended to general generated Jacobian equations in \cite{Je}.  As a byproduct of our local convexity theory in \cite{LoT2}, these results also extend to $C^2$ cost and generating functions.

\vspace{2mm}

\noindent {\bf 4.2 Other convexity results}. We prove here a version of Lemma 2.2 which does not use duality and corresponds in some sense to Lemma 2.3. We use the same notation as in Lemma 2.2.

 \begin{lem}
Assume $g$ satisfies A1,A2,A1* and A3w  and suppose $u\in C^0(\Om)$ is $g$-convex in $\Om$.  Assume also:
\vspace{2mm}

(i) $\Omega$ is $g$-convex with respect to some $ \tilde y\in \Sigma_0$,  $ \tilde z = g^*(x_0, \tilde y, u_0)$ for some $x_0$ in $\Om$; 
\vspace{2mm}

(ii) $Tu(x_0)$  is sub $g^*$-convex with respect to $\tilde g =g(\cdot, \tilde y, \tilde z)$ in $\Omega$.
\vspace{2mm}

\noindent Then $\tilde y \in Tu(x_0)$.

\end{lem}

To prove Lemma 4.3 we first choose two extreme points $p_0$ and $p_1$ in $\partial u(x_0)$ such that $p_\theta = \theta p_0 + (1-\theta) p_1 =g_x(x_0, \tilde y,  \tilde z)$ for some  $\theta\in (0,1)$. Setting $y_\theta , z_\theta = Y, Z(x_0, u_0, p_\theta)$, and $g_\theta = g(\cdot, y_\theta,z_\theta)$  for any $\theta \in [0,1]$, it follows that
for either $i = 0$ or $1$, $\tilde  h = g_i-\tilde g > 0$ at  some point on the $g$-segment $I_g$, with respect to $\tilde y, \tilde z$, joining $x_0$ and a point $x\in \Omega$, provided $D_\eta g_0(x_0) \ne D_\eta g_1 (x_0)$, where $\eta_j = E^{i,j}(x_0, \tilde y,\tilde z)(q_i - q_{0i})$, 
$q = Q(x,\tilde y, \tilde z), q_0 = Q(x_0, \tilde y,\tilde z)$. By decreasing $\tilde z$ we can then apply the argument of Lemma 2.3 in the domain $\Omega$ to obtain a contradiction if $ \max\{g_0,g_1\}(x) < \tilde g(x)$. On the other hand, if $D_\eta g_0(x_0) = D_\eta g_1 (x_0)$, then the function $f$, given by $f(t)= \tilde h(x_t)$, 
 satisfies  $f(0) = f^\prime (0) = 0$ for both $i=0$ and $1$ so that if also $f\le0$ on $[0,1]$, $f(1)< 0$, we also obtain a contradiction with the differential inequality, \eqref{f diff ineq}. Alternatively, we can approximate $x$ to reduce to the case 
 $D_\eta g_0(x_0) \ne D_\eta g_1 (x_0)$, as in \cite{T5}. Consequently $\tilde g \le \max\{g_0,g_1\}$ on $\Om$ whence $\tilde y\in Tu(x_0)$ as asserted. 

Note that the  domain $\Om$ in Lemma 4.3 can be made arbitrary by replacing it by its $Q(\cdot, \tilde y, \tilde z)$ convex hull in condition (ii). By extending to all $\tilde y \in \Sigma _0$, we can then obtain another version of Lemma 2.2, though under stricter sub convexity hypotheses.
 \vspace{0.2cm}

We also have another version of Lemma 2.3 from our earlier drafts, which essentially follows from the proof of Lemma 2.1. 

 \begin{lem}
Assume $g$ satisfies A1,A2,A1* and A3w  and suppose $u\in C^0(\Om)$ is $g$-convex in $\Om$ and $g_0 = g(\cdot, y_0, z_0)$ is a $g$-affine function on $\Om$.  Assume also:
\vspace{2mm}

(i) $\Omega$ is $g$-convex with respect to  $y_0$ and all $z\in g^*( \cdot, y_0, \sup\{ u,g_0\})(\Om)$;
\vspace{2mm}

(ii) $Tu(\Omega)\cup \{y_0\}$  is sub $g^*$-convex with respect to 
$\sup\{u,g_0\}$ in $\Omega$.
\vspace{2mm}

\noindent Then the section $S = S(u,g_0)$ is $g$-convex with respect to $(y_0, z_0)$, while if $g_0$ is a $g$-support to $u$, the contact set $S_0 = S_0(u,g_0)$ is $g$-convex with respect  to $(y_0, z_0)$. 

\end{lem}

Note that when we extend Lemma 4.4 to general domains, when $g_0$ is an arbitrary $g$-support to $u$, we end up with the same sub convexity hypotheses as the corresponding extension above of Lemma 2.3.   

Finally we also indicate that by modification of the proof of Lemma 4.3, we may obtain a refinement of the local version of the Loeper maximum principle \eqref{Loeper} under condition A3, corresponding to those used to show $C^{1,\alpha}$ regularity in \cite{Loe, GuT2015, Je}. This is  also extended to $g \in C^2$, through the geometric approach in \cite{LoT2}.

\begin{lem}
Assume $g$ satisfies A1,A2,A1* and A3 and  $(x_0, u_0, [p_0,p_1]) \subset \mathcal U^\prime \subset\subset \mathcal U$, where $[p_0,p_1]$ denotes the straight line joining $p_0$ and $ p_1$ in $\RR^n$. Then there exist small positive constants,  $\epsilon_0$ and $\gamma_0$, depending on $g$ and $\mathcal U^\prime$, such that

\begin{equation} \label{4.4}
g(x, y_\theta, z_\theta)  \le \max\{ g_0(x), g_1(x)\} - \gamma_0 [\theta(1-\theta) |p_1-p_0| |x-x_0|]^2
\end{equation}

\noindent for all $|x-x_0| \le \epsilon_0$ and $\theta\in [0,1]$.
\end{lem}

To prove \eqref{4.4}, we replace $\tilde g = g_\theta$ in the proof of Lemma 4.3 by the function $\tilde g_\gamma$ given by

\begin{equation} \label{4.5}
\tilde g_\gamma  = g_\theta + \gamma [\theta(1-\theta) |p_1-p_0| |x-x_0|]^2
\end{equation}

\noindent for sufficiently small $\gamma > 0$  and distance $|x-x_0|$, using more explicitly the differential inequalities (2.9) and (2.10) in \cite{T5}, in conjunction with  condition A3 in $\mathcal U^\prime$, to extend \eqref{f diff ineq} to $\tilde g_\gamma$.


\vspace{1cm}



\frenchspacing
\begin{thebibliography}{99}

\bibitem{GT} 
\newblock D. Gilbarg and N.S. Trudinger,
\newblock Elliptic Partial Differential Equations of Second Order, 
\newblock Springer, 1983.



\bibitem{GK}
\newblock N. Guillen and J. Kitagawa,
\newblock \emph{Pointwise inequalities in geometric optics and other generated Jacobian equations},
\newblock  Comm.\ Pure Appl.\ Math. \textbf{70} (2017), 1146-1220.
.

\bibitem{GuT2015} 
\newblock C.E. Guti\'errez and  F. Tournier, 
\newblock \emph{Regularity for the near field parallel refractor and reflector problems},  \textbf{45}, 917-949 (2015)
\newblock Calc.\ Var. \ Partial Differ. \ Equ. \textbf{45},(2015), 917-949.

\bibitem{Je} 
\newblock S. Jeong,
\newblock \emph {Local Holder regularity of solutions to generated Jacobian equations},
\newblock preprint, arXiv, 2020.

\bibitem{JT1} 
\newblock F. Jiang and N. S. Trudinger,
\newblock \emph {On Pogorelev estimates in optimal transportation and geometric optics},
\newblock Bull.\ Math.\ Sci. \textbf{4} (2014), 407-431.

\bibitem{JT2} 
\newblock F. Jiang and N. S. Trudinger,
\newblock \emph {On the second boundary value problem for Monge Amp\`ere type Equations and geometric optics},
\newblock Arch.\ Rational Mech.\ Anal. \textbf{229} (2018), 547-567.

\bibitem{JTY} 
\newblock F. Jiang, N. S. Trudinger and X.-P. Yang,
\newblock \emph {On the Dirichlet Problem for Monge-Amp\`ere type equations},
\newblock Calc.\  Var.\  Partial Differ.\  Equ. \textbf{49} (2014), 1223-1236.


\bibitem{KW}
\newblock   A. Karakhanyan and X.-J. Wang,
\newblock    \emph{On the reflector shape design},
\newblock   J. \ Diff. \ Geom.  \textbf{84} (2010), 561--610.

\bibitem{KM}
\newblock Y.-H. Kim and R.J. McCann,
\newblock \emph{Continuity, curvature and the general covariance of optimal transportation},
\newblock  J. \ Eur. \ Math. \ Soc. \textbf{12} (2010), 1009--1040.
  
 \bibitem{LT1}
\newblock      J.-K. Liu and N.S. Trudinger,
\newblock     \emph{On Pogorelov estimates for Monge-Amp\`ere type equations},
\newblock      Discrete Contin.\ Dyn.\ Syst.\ Ser. A. \textbf{28}, (2010), 1121--1135.

 \bibitem{LT2}
\newblock      J.-K Liu and N.S Trudinger,
\newblock    \emph{On classical solutions of near field reflection problems},
\newblock    Discrete Contin. Dyn. Syst. \textbf{36},  (2016), 895-916.
    
 \bibitem{LTW} 
\newblock J.-K. Liu, N.S Trudinger and X.J. Wang, 
\newblock     \emph{Interior $C^{2,\alpha}$ regularity for potential functions in optimal transportation}.
\newblock      Comm.\ Partial Differential Equations \textbf{35} (2010), 165--184.

\bibitem{Loe}
\newblock     G. Loeper,
\newblock    \emph{On the regularity of solutions of optimal transportation problems},
\newblock   Acta Math. \textbf{202} (2009), 241--283.

 \bibitem{LoT1}
\newblock      G. Loeper and N.S Trudinger,
\newblock     \emph{Weak formulation of the MTW condition and convexity properties of potentials},
\newblock   preprint, arXiv, 2020, Methods Appl. Anal. (to appear).

 \bibitem{LoT2}
\newblock      G. Loeper and N.S Trudinger,
\newblock     \emph{On the convexity theory of generating functions},
\newblock  preprint (2020).

\bibitem{MTW} 
\newblock X.-N. Ma, N.S.Trudinger and X.-J.Wang,
\newblock   \emph{Regularity of potential functions of the optimal transportation problem},
\newblock   Arch.\ Rat.\ Mech.\ Anal.\  \textbf{177} (2005), 151-183.

\bibitem{Ra} 
\newblock C. Rankin, 
\newblock \emph{Distinct solutions to GPJE cannot intersect},
\newblock Bull.\ Aust.\ Math.\ Soc.\ \textbf{102} (2020),  462-470.


\bibitem{T1}  
\newblock N.S.Trudinger,
\newblock  \emph {Recent  developments in elliptic partial differential equations of  Monge-Amp\`ere type},
\newblock Proc.\ Int.\ Cong.\ Math.,Madrid, \textbf{3} (2006), 291-302.

 \bibitem{T4}
 \newblock   N.S.Trudinger, 
 \newblock  \emph {On generated prescribed Jacobian equations},
 \newblock  Oberwolfach Reports \textbf{38} (2011), 32-36.
     
 \bibitem{T5}
 \newblock  N.S.Trudinger,
 \newblock  \emph {The local theory of prescribed Jacobian equations}, 
 \newblock Discrete Contin.\ Dyn. \ Syst. \textbf{34} (2014), 1663-1681.
   
 
 \bibitem{TW1}
\newblock     N.S.Trudinger and X.-J. Wang,
\newblock     \emph{The Monge-Amp\`ere equation and its geometric applications},
\newblock     in "Handbook of Geometric Analysis", International Press (2008), 467-524.

\bibitem{TW2}
\newblock    N.S.Trudinger and X.-J. Wang,
\newblock     \emph{On convexity notions in  optimal transportation},
\newblock     preprint  (2008).

\bibitem{TW3}
\newblock    N.S.Trudinger and X.-J. Wang,
\newblock   \emph{On the second boundary value problem for Monge-Amp\`{e}re type equations and optimal transportation},
\newblock   Ann.\ Scuola Norm.\ Sup.\ Pisa Cl.\ Sci.  \textbf{VIII}, (2009), 143--174.

\bibitem{TW4}   
\newblock N.S.Trudinger and X.-J.Wang,
\newblock    \emph{On strict convexity and continuous differentiability of potential functions in optimal transportation},
\newblock   Arch.\ Rat.\ Mech.\. Anal. \textbf{192} (2009), 403--418.

\bibitem {V2} 
\newblock C.Villani,
 \newblock Optimal Transportation, Old and New,
          Springer, 2008

\end{thebibliography}
\end{document}